\documentclass{article}        % Your input file must contain these two lines 
\usepackage{amsmath}
\usepackage{amssymb}

\begin{document}               % plus the \end{document} command at the end.

\newcommand{\forces}{\Vdash}
\newcommand{\tv}[1]{[\![#1]\!]}

\title{A Cheerful Introduction to Forcing and the Continuum
  Hypothesis}\author{Kenny Easwaran}
  \date{}\maketitle{}

This text grew out of a presentation I made in the Berkeley math
department's ``Many Cheerful Facts'' series on November 16, 2005.  These
talks are intended to be accessible to most math graduate students, to
help them understand important concepts and connections that exist in
areas other than the one(s) in which they are working.  In general,
these talks presuppose familiarity with concepts like groups,
topological spaces, categories and the like.  But in practice, people
giving talks in more ``core'' areas of the Berkeley department (like
number theory and algebraic geometry) have been known to presuppose even
more material.  Since my talk was on a subject very different from any
of this, I tried to presuppose as little as possible about set theory
(mainly just a familiarity with the notation of membership ``$\in$'',
subsethood ``$\subseteq$'', the union and intersection operations, and
the central fact that there is no one-to-one correspondence
(``bijection'') between the elements of a set and its powerset).  In the
talk I took advantage of the basic algebraic knowledge of my audience to
ease the discussion of boolean algebras, but since I would like this
written version to be useful to philosophers as well as mathematicians,
I have included two versions of this discussion.

The material I cover here is definitely too technical to cover in full
detail to an audience without a background in set theory.  However, the
important technicalities are only relevant for proving certain lemmas
that of little interest for understanding the results.  I have
explicitly mentioned any point where I have omitted such a proof - they
can all be found in standard references like \cite{bell}, \cite{jech},
and \cite{kunen}.  At many other points, the technicalities are not so
complicated as to warrant leaving them out entirely, but I felt that
they would interrupt the explanation.  In these cases, I have made
extensive use of footnotes to fill in details of various arguments.
Ideally, the reader should be able to ignore all the footnotes, and only
consult them when she has a particular interest in the details of
particular segments of the argument.  This may make the structure of the
document somewhat confusing, but I hope that it makes it easier for
people with different backgrounds to read it.  Every exposition of
forcing that I have seen either presumes a fair bit of familiarity with
set theory, or omits all the discussion of how the method works.  I hope
to fill this gap.

\section{Boolean algebras and partial orders}\label{BAs}

In this section, I will explain what you need to know about boolean
algebras (including their ordering, and the operations
$\wedge,\vee,\neg$), ultrafilters, and complete boolean algebras in
order to understand the method of forcing.  If you are already familiar
with all those terms, just skip to subsection \ref{pos} to learn a few
facts about separative, atomless partial orders.

\subsection{Boolean algebras}

\subsubsection{The basics}

A partially ordered set is a set of elements with a relation $x\leq y$
that is reflexive, antisymmetric, and transitive.  That is, $x\leq x$,
and if $x\leq y$ and $y\leq x$ then $x=y$, and if $x\leq y$ and $y\leq
z$ then $x\leq z$.  (A linear ordering, or total ordering, is a partial
ordering with the additional requirement of trichotomy, that for any $x$
and $y$ either $x\leq y$ or $y\leq x$ or $x=y$.)  A boolean algebra is a
particular sort of partially ordered set (or ``partial order'', or
``poset'', for short).  A boolean algebra must have a maximal and a
minimal element, denoted by ``$1$'' and ``$0$'' respectively.  (That is,
for any $x$, $1\geq x\geq 0$.)  In addition, every pair of elements $x$
and $y$ must have a greatest lower bound $x\wedge y$ and least upper
bound $x\vee y$.  (That is, $x\wedge y\leq x$ and $x\wedge y\leq y$, and
if $z\leq x$ and $z\leq y$, then $z\leq x\wedge y$; and dually for
$x\vee y$.)

Given these requirements, we can see that $x\vee y=y\vee x$ and $x\wedge
y=y\wedge x$, (commutativity) $x\vee(y\vee z)=(x\vee y)\vee z$ and
$x\wedge(y\wedge z)=(x\wedge y)\wedge z$, (associativity) $x\vee 1=1$,
$x\wedge 1=x$, $x\vee 0=x$, $x\wedge 0=0$, and finally $x\leq y$ iff
$x\wedge y=x$.

We also require that for every $x$, there is an element $\neg x$ such
that $x\vee\neg x=1$ and $x\wedge\neg x=0$.  Finally, we require that
$x\vee(y\wedge z)=(x\vee y)\wedge(x\vee z)$ and $x\wedge(y\vee
z)=(x\wedge y)\vee(x\wedge z)$ (distributivity).  Using these facts, we
can show that $\neg x$ is unique.\footnote{If $x\vee y=1$ and $x\wedge
  y=0$, then $(x\vee \neg x)\wedge y=1\wedge y=y$, so $y=(x\wedge
  y)\vee(\neg x\wedge y)=0\vee(\neg x\wedge y)=\neg x\wedge y$, so
  $y\leq\neg x$. Similarly, $(x\wedge\neg x)\vee y)=0\vee y=y$, so
  $y=(x\vee y)\wedge(\neg x\vee y)=1\wedge(\neg x\vee y)=\neg x\vee y$.
  Thus, $\neg x\wedge y=\neg x\wedge(\neg x\vee y)=(\neg x\wedge\neg
  x)\vee(\neg x\wedge y)=\neg x\vee z$ where $z\leq \neg x$, so $\neg
  x\vee z=\neg x$.  Thus, $\neg x\wedge y=\neg x$, so $\neg x\leq y$,
  and by antisymmetry, we see that $y=\neg x$.}  It is clear that
$\neg\neg x=x$, and it is not hard to show that $\neg(x\vee
y)=(\neg x)\wedge(\neg y)$, and $\neg(x\wedge y)=(\neg x)\vee(\neg y)$
(DeMorgan's Laws).\footnote{$((\neg x)\wedge(\neg y))\vee(x\vee
  y)=((\neg x)\vee(x\vee y))\wedge((\neg y)\vee(x\vee y))=(1\vee
  y)\wedge(1\vee x)=1$ and $((\neg x)\wedge(\neg y))\wedge(x\vee
  y)=(((\neg x)\wedge(\neg y))\wedge x)\vee((\neg x)\wedge(\neg y)\wedge
  y)=(0\wedge \neg y)\vee(0\wedge\neg x)=0$.}  Thus, we can see that
$x\leq y$ iff $\neg y\leq \neg x$.\footnote{It is easy to check that if
  we define $x\cdot y=x\wedge y$ and $x+y=(x\wedge\neg
  y)\vee(y\wedge\neg x)$, then a boolean algebra is a commutative ring
  with identity such that $\forall x(x\cdot x=x)$.  One can also verify
  that in any such ring, if we define $x\leq y$ to mean $x\cdot y=x$,
  then the resulting structure is a boolean algebra.  It turns out that
  $x\wedge y=x\cdot y$ and $x\vee y=x+y+x\cdot y$.}

\subsubsection{Ideals and filters, and examples}

An ideal is a nonempty, proper subset of a boolean algebra that is
closed under $\vee$ and closed downwards under the order.  That is, if
$x,y\in I$ and $z\leq x$ then $z\in I$ and $x\vee y\in I$.\footnote{It
  is easy to verify that this definition of an ideal is equivalent to
  the ring-theoretic definition.}  A filter is just the dual - if
$x,y\in F$ and $z\geq x$ then $z\in F$ and $x\wedge y\in F$.  It is
straightforward to check that $0$ is in every ideal but no filter, and
$1$ is in every filter but no ideal.  Thus, if $x$ is in an ideal
(filter), then $\neg x$ is not in the ideal (filter).  If the converse
also holds (that is, that if $\neg x\not\in I$ implies $x\in I$, or
similarly for a filter $F$), then the set is said to be a prime ideal,
or an ultrafilter.\footnote{Note that $x\cdot\neg x=x\wedge\neg x=0$, so
  since $0\in I$ for any ideal, a prime ideal must satisfy the above
  requirement.  The converse is also true.  It's also easy to see that
  any such ideal is maximal, because if $x\not\in I$, then $\neg x\in
  I$, so adding $x$ to $I$ would also put $x+\neg x=1$ into the ideal.
  The fact that every prime ideal is maximal can also be recognized by
  noting that every element other than $1$ is a zero divisor.}

For any ideal $I$, the set $\{x|\neg x\in I\}$ is a filter, and vice
versa.  If $I$ is a prime ideal, then this set is an ultrafilter, and
vice versa.

The natural example of a boolean algebra is a formal language for logic,
where the elements are equivalence classes of formulas rather than
formulas themselves.  On any interpretation of such a language, the set
of true sentences forms an ultrafilter, and the set of false ones forms
a prime ideal.  Given a collection of interpretations, the set of
sentences true on all of them forms a filter, and the set of sentences
false on all of them forms an ideal.

Dually, if we think of the boolean algebra as a set of possible ``truth
values'' to assign to sentences, and then specify some of them to count
as ``true'', then the specified set will be an ultrafilter.

If we have one boolean algebra $B$ and the function $f$ maps it into
another boolean algebra $C$ in such a way that the operations
$(\wedge,\vee,\neg)$ are all preserved by the mapping (that is,
$f(x\wedge y)=f(x)\wedge f(y)$, and so on), then the set of elements
mapped to $1$ forms a filter and the set of elements mapped to $0$ forms
an ideal.\footnote{This is just the standard notion of a ring
  homomorphism, so obviously the kernel is an ideal.}  (Any further
properties I mention of ideals and $0$ naturally generalize to
properties of filters and $1$, with perhaps a few other obvious
changes.)  Conversely, if $B$ is a boolean algebra and $I$ is an ideal,
then there is a boolean algebra $B/I$ (``$B$ modulo $I$'', or ``the
quotient of $B$ by $I$'') with a natural map $f$ from $B$ to $B/I$ such
that $f(x)=0$ iff $x\in I$.  $B/I$ can be thought of as the set of
equivalence classes of elements of $B$, where two elements are said to
be equivalent just in case $x\wedge\neg y$ and $y\wedge\neg x$ are both
in $I$.  In addition, if $g\colon B\rightarrow C$ and $g(x)=0$ iff $x\in
I$, then the set of elements in the range of $g$ is isomorphic to
$B/I$.

Now consider any set $S$, and its powerset $\mathcal{P}(S)$.  If we let
$x\leq y$ mean $x\subseteq y$, then it is straightforward to see that
this is a boolean algebra, with $x\wedge y=x\cap y$, $x\vee y=x\cup y$,
and $\neg x=S\setminus x$.  This algebra is said to be a
\emph{complete} boolean algebra, because \emph{every} set of elements
has a least upper bound (their union) and a greatest lower bound (their
intersection), not just the pairs, as required for an ordinary boolean
algebra.\footnote{This is the same notion of completeness that the reals
  have - the rationals are such that any pair of elements has a least
  upper bound, but the reals have the additional property that any
  \emph{bounded}, not just finite, set of elements has a least upper
  bound.  In a boolean algebra, every set of elements is bounded by $1$,
  so to be complete, \emph{every} set must have a least upper bound.}
In this algebra, the singleton sets $\{a\}$, where $a\in S$, are said to
be atoms - that is, if $x\leq\{a\}$, then either $x=0$ or
$x=\{a\}$. Every element has an atom below it, so the algebra is said to
be atomic. A boolean algebra without atoms is said to be atomless.

Let $S$ be an infinite set, and then let $I$ be the collection of its
finite subsets.  Then $I$ is an ideal, because the union of two finite
sets is finite, and any subset of a finite set is finite.  The quotient
$\mathcal{P}(S)/I$ will then be an atomless boolean algebra.  This is
because if $x$ is non-zero, then it must be an equivalence class
containing some infinite subset of $S$.  (Every finite subset is
equivalent to $0$ in the quotient.)  If $y$ and $z$ are any two disjoint
infinite subsets of $x$, then neither is equivalent to either $x$ or
$0$, because they differ by infinitely many elements from each.  Thus,
$x$ is not an atom.  This boolean algebra is not complete, because if we
let $x_1,x_2,\dots$ be countably many disjoint infinite subsets of $S$,
then they have no least upper bound.  Any mutual upper bound $x$ must
contain all but finitely many elements of each of these sets.  But if we
now let $e_1\in x_1\cap x, e_2\in x_2\cap x,\dots$, then
$x\setminus\{e_1,e_2,\dots\}$ is also a mutual upper bound of the $x_i$,
but it is distinct from $x$ (since they differ on infinitely many
elements) and is below it in the ordering (since it is a subset of $x$).
Thus, there is no least upper bound, so the algebra is not complete.

\subsection{Partial orders}\label{pos}

When dealing with a boolean algebra $B$, the partially ordered set that
we will consider is normally $B\setminus\{0\}$.  In this partial order
or in any other, we define $x\perp y$ iff $x$ and $y$ have no common
lower bound - that is, there is no $z$ such that $z\leq x$ and $z\leq
y$.  An atomless boolean algebra gives a partial order such that
$\forall x\exists y(y<x)$, and this is how we define an atomless partial
order in general.

One other important fact about boolean algebras is that if $x\not\geq y$
(so that $x\wedge y\not=y$) then,
since $y=(y\wedge x)\vee(y\wedge\neg x)$, we see that $y\wedge\neg
x\not=0$.  If we let $z=y\wedge\neg x$, then we see that if $x\not\geq y$ then
$\exists z(z<y\text{ and }z\perp x)$.  This property holds for all
boolean algebras, and if it holds in a general partial order, we call
such a partial ordering ``separative''.  The idea is that if we think of
think of each element of the partial ordering as a piece of information
about how some world might be (it could be a proposition about the
world, or a set of possible worlds, or could somehow specify such a
set), then unless $y$ entails $x$ (in which case $y\leq x$), there is
some piece of information $z$ extending $y$ that is incompatible with
$x$.  Any piece of information can be extended to two incompatible
pieces of information.  Moving downwards in the partial ordering always
corresponds to getting more information about this possible world, since
there are always other possibilities that are being ruled out.

In an atomless partial order $P$, we can say that a set $D$ is dense iff
$$\forall p\exists q(q\in D\text{ and }q\leq p),$$
so that every element of the ordering can be ``refined'', or
``extended'', to get a further element of the set.  We will later be
concerned with ultrafilters that intersect various dense subsets of our
boolean algebra.  However, it is useful to note that if $P$ is a partial
order that is a dense subset of some boolean algebra $B$, and $D$ is a
dense subset of $P$, then $D$ is also a dense subset of $B$.  In
addition, if $D$ is a dense subset of $B$, then $D'=\{x\in P|\exists
y(x\leq y\text{ and }y\in D\}$ is a dense subset of $P$.  Thus, $P$ and
$B$ are in some senses interchangeable when considering their
collections of dense subsets.  In particular, any filter intersecting
every dense subset of one corresponds to a filter intersecting every
dense subset of the other.  Such filters play the central role in
sections \ref{names} and \ref{forcing}.

The important fact about separative, atomless partial orders is that if
$P$ is such an ordering, then there is a \emph{unique} complete boolean
algebra $B$ such that $P$ is isomorphic (as an ordering) to a dense
subset of $B$.\footnote{The construction involves considering the
  collection of regular open subsets of $P$ under the order topology.
  The technical details are unimportant here.}

\section{The set-theoretic world-view}\label{settheory}

In set theory, we pretend that there are no objects other than
sets.\footnote{Some set theorists actually believe this, at least about
  mathematical objects - they think natural numbers, vector spaces,
  schemes, etc. are all just sets of certain sorts.  Actually adopting
  this attitude is not necessary here - we just need to restrict
  attention to the sets alone, if there is anything else.}  This means that
there is the empty set, $\emptyset$; the set containing that,
$\{\emptyset\}$; the set containing both of those,
$\{\emptyset,\{\emptyset\}\}$; and the like.  Set theorists identify
the empty set with the number $0$, and the other two sets mentioned
above with $1$ and $2$ respectively.  In general, we can let
$n=\{0,\dots,n-1\}$, so that the natural numbers have a canonical
representation as sets.

In addition to each of the natural numbers, there is also a set
containing all of them.  Mathematicians standardly call this set
$\mathbb{N}$, but set theorists call it $\omega$ for technical and
historical reasons.  Once the natural numbers have been identified, the
integers can be represented as ordered pairs of naturals, rationals as
ordered pairs of integers, and reals as Dedekind cuts of
rationals.\footnote{Ordered pairs can be represented as sets as follows:
  $(x,y)=\{\{x\},\{x,y\}\}$.  With this identification, we can see that
  (assuming $x$ and $y$ are distinct) each ordered pair will have one
  singleton and one unordered pair as elements.  Whichever element is in
  the singleton is thus represented as being ``the first'' of the
  two. It's easy to check that with this identification, $(x,y)=(z,w)$
  iff $x=z$ and $y=w$, even if $x=y$ (since then $z$ will have to equal
  $w$).

With this identification, we then represent integers as equivalence
classes of ordered pairs where $(n_1,n_2)=(m_1,m_2)$ iff
$n_1+m_2=n_2+m_1$ - the equivalence class of the pair $(n_1,n_2)$ is
then seen as coding the integer $n_1-n_2$.  The coding of rationals as
pairs of integers is more familiar - $p/q$ is coded by the equivalence
class of the pair $(p,q)$, and two pairs are equivalent just in case
$p_1q_2=q_1p_2$.  A Dedekind cut is a set $r$ of rationals such that if
$q<q'$ and $q'\in r$, then $q\in r$.  Each such cut codes the real that
would be the least upper bound of the set.} Other mathematical entities
can be represented using similar means, so there really is no
restriction to saying that sets are the only things that exist.

In the universe of sets, every set has a powerset, which is the set of
all subsets of the first set.  The powerset of $X$ is symbolized as
$\mathcal{P}(X)$.  If two sets can be put in bijection with one
another,\footnote{A bijection is just a one-to-one and onto function.  A
  function $f\colon A\rightarrow B$ is a bijection iff for all $x\in B$
  there is a unique $y\in A$ such that $f(y)=x$.} we will say they have
the same cardinality ($|x|=|y|$).  Some well-known results of Cantor
show that no set can be put in bijection with its powerset, and that the
powerset of $\omega$ can be put in bijection with the set of
reals.\footnote{All of this is covered in many places.  To see that
  there is no bijection between $S$ and $\mathcal{P}(S)$, let that
  bijection be $f$, and then consider the set $\{x\in S\colon x\not\in
  f(x)\}$, which can't be in the range of the supposed bijection $f$.
  To see that the reals and the powerset of the naturals are
  equinumerous, pair each set of naturals with the real number whose
  $n$th bit in binary is 1 iff $n$ is in the set.  Some care is needed
  for reals that can be represented in two different ways in binary,
  like $.10000\dots$ and $.01111\dots$.} Therefore, I will use
$\mathbb{R}$ to refer to $\mathcal{P}(\omega)$, since the relevant issue
here is just cardinality.\footnote{Set theorists often let
  ``$\mathbb{R}$'' denote the set of functions from $\omega$ to
  $\omega$, or from $\omega$ to $2$ for the same reason, as well as
  further topological similarities.  At any rate, it's generally worth
  double-schecking what a set theorist means when she talks about a real
  number.}

Because every set has a powerset, and we know that a powerset has a
strictly larger cardinality than its corresponding set, we know there
are infinitely many infinite cardinalities.  Through use of the axiom of
choice, we see that any two infinite cardinalities can be compared, and
they can in fact be well-ordered,\footnote{That is, in any non-empty
  set of cardinals, there is a unique least element.} so we can call the
infinitely cardinalities in order $\aleph_0,\aleph_1,\aleph_2,\dots$, where
$\aleph_0$ is the cardinality of the naturals.\footnote{Technically,
  we'll need the notion of an ordinal to be able to number all the
  infinite cardinalities.  In the context of the Axiom of Choice,
  cardinals are generally identified with certain ordinals.  Thus, they
  can be well-ordered, so it really does make sense to talk about ``the
  least'' uncountable cardinal. The technical details of this can be
  found in any book on set theory, for instance Kunen, ch. 1.}  Since we
know the reals have a larger cardinality then $\aleph_0$, and can't
explicitly construct any sets of intermediate size, we might conjecture
that $|\mathbb{R}|=\aleph_1$.  In fact, this is exactly what Cantor did,
and this hypothesis is called the Continuum Hypothesis (CH).  The
Generalized Continuum Hypothesis (GCH) states that for every set $x$
with cardinality $\kappa$, the powerset of $x$, $\mathcal{P}(x)$, has
cardinality $\kappa^+$, the least cardinal greater than $\kappa$.  What
I will prove below is that if the known axioms of set theory (ZFC) are
consistent, then so is ZFC+$|\mathbb{R}|=\aleph_\alpha$ for just about
every $\alpha$, so we cannot prove CH from ZFC.\footnote{Clearly, we
  can't have $\alpha=0$, by Cantor's theorem.  In addition, there are
  certain other values for $\alpha$, like $\omega$, that are ruled out
  by a result called K\"{o}nig's Theorem, even though both larger and
  smaller values are possible.}  The proof can easily be generalized to
show that GCH is violated at any particular cardinality, not just at
$\aleph_0$.

\section{Models of set theory}\label{models}

Set theorists (like all logicians) like to formalize their vocabulary.
In this case, the statements of interest can all be built up from
sentences of the form $x=y$ and $x\in y$ by use of logical connectives
(and, or, not, iff, etc.) and quantifiers ($\forall x$ and $\exists
x$).  A typical example is the sentence
$$\forall x\forall y(\exists z(\forall w(w\in z\leftrightarrow
(w=x\text{ or }w=y))))$$
This sentence is one of the axioms of ZFC that are taken to be the facts
that we know about the universe of sets.  It states that for any sets
$x$ and $y$, the set $\{x,y\}$ exists as well.\footnote{If you want the
  gory details, the complete list is as follows.  There are other
  versions of some of these, and slightly different overall sets.  Any
  good book on set theory will have them listed somewhere in an early
  chapter, if not on the first page. Most of the details aren't too
  important here, though they would be for the detailed proofs of the
  results I mention.\begin{enumerate}
\item Extensionality - $\forall x\forall y (x=y\leftrightarrow\forall
  z(z\in x\leftrightarrow z\in y))$
\item Pairing - $\forall x\forall y(\exists z(\forall w(w\in z\leftrightarrow
(w=x\text{ or }w=y))))$ - $z$ is called $\{x,y\}$
\item Separation schema - $\forall x(\exists y(\forall z(z\in y\leftrightarrow
  (z\in x\text{ and }\phi(z)))))$ - $y$ is called $\{z\in y|\phi(z)\}$
\item Union - $\forall x(\exists y(\forall z(z\in
  y\leftrightarrow\exists w(z\in w\text{ and } w\in x))))$ - $y$ is
  called $\bigcup x$
\item Powerset - $\forall x(\exists y(\forall z(z\in
  y\leftrightarrow\forall w(w\in z\rightarrow w\in x))))$ - $y$ is
  called $\mathcal{P}(x)$
\item Replacement schema - $\forall x(\exists y(\forall z(z\in y\leftrightarrow
  \exists w(w\in x\text{ and }\phi(w,z)))))$, for any formula
  $\phi(w,z)$ such that $\forall w\exists ! z(\phi(w,z)$ - that is, such
  that for any $w$ there is a unique $z$ satisfying $\phi(w,z)$, so that
  $\phi$ represents a function.
\item Infinity - $\exists y(\emptyset\in y\text{ and }\forall x(x\in
  y\rightarrow (x\cup\{x\})\in y))$ - the smallest such $y$ is called $\omega$
\item Foundation - $\forall x(\exists y(y\in x\text{ and }\forall z(z\in
  x\rightarrow z\not\in y)))$
\item Choice - $\forall x[(\forall y(y\in x\rightarrow\exists z(z\in
  y)))\rightarrow\exists f((f\colon x\rightarrow \bigcup x)\text{ and
  }\forall z(z\in x\rightarrow (f(z)\in z)))]$
\end{enumerate}
  In Separation and Replacement, $\phi$ is any arbitrary formula with
  just the stated variables free that can be written in this language,
  so those two are actually infinite sets (schemas) of axioms rather
  than individual axioms. In Choice, the notation $f\colon x\rightarrow
  y$ is an abbreviation for the (very long) sentence saying that $f$ is
  a set of ordered pairs whose first elements are all in $x$, second
  elements are all in $y$, and such that each element of $x$ is the
  first element of exactly one of the pairs.} Basically anything of
set-theoretic interest can be phrased in this sort of language.  Given
any sentence $\phi$ in this notation, and any transitive set $M$ (a
transitive set is one that contains all elements of its elements, so if
$x\in y\in M$ then $x\in M$), we can find a related sentence known as
$\phi^M$ where all quantifiers are restricted to $M$.\footnote{The
  restriction works as follows - $\forall x \phi(x)$ becomes $\forall
  x(x\in M\rightarrow \phi(x))$ and $\exists x\phi(x)$ becomes $\exists
  x(x\in M\text{ and }\phi(x))$.  This works just as one would expect.}
If the above sentence is taken as $\phi$, asserting that for any sets
$x$ and $y$, $\{x,y\}$ exists, then the restriction $\phi^M$ asserts
that if $x$ and $y$ are in $M$, then $\{x,y\}\in M$ as
well.\footnote{Actually, it states that there is some element of $M$
  whose only elements \emph{in $M$} are $x$ and $y$ - but since $M$ is
  transitive, all elements of $z\in M$ must be elements of $M$ as well,
  and since the restricted axiom states that the only elements of $z$ in
  $M$ are $x$ and $y$, we see that $z$ must in fact be $\{x,y\}$.}  Once
we have defined $\phi^M$, we say that $M\models\phi$ (in words, ``$M$
satisfies $\phi$'', or ``$M$ is a model of $\phi$''\footnote{The notion
  of model used here is much like the notion used in model theory, but
  is not quite the same.  In model theory, instead of just fixing a set
  and letting the $\in$ relation in the model be the actual $\in$
  relation, one must specify a set of ordered pairs to stand for the
  $\in$ relation in addition to specifying the domain of quantification.
  This allows a unified treatment of models of any sort of theory, not
  just theories phrased in the language of set theory.  In addition, it
  allows for the construction of models with ``non-standard'' natural
  numbers and the like.  But for our purposes, we must stick with models
  where $\in$ represents the actual $\in$ relation, and enough of the
  relevant results from model theory will carry through.  This
  restriction to the actual $\in$ relationship is why we only consider
  restrictions to transitive sets $M$.}) just in case $\phi^M$
is true.  The important point about this relationship is that
G\"{o}del's Completeness Theorem guarantees that for any set $T$ of
formulas in the language of set theory, this set is consistent iff there
is some set $M$ such that $M\models T$.  In particular, if we assume
that ZFC is consistent, then there is some $M$ such that
$M\models\mathrm{ZFC}$.\footnote{Technically, we can't guarantee that
  this set $M$ is transitive. However, using some set-theoretic trickery
  (the Reflection Theorem schema), we can guarantee that if $T$ extends
  ZFC, then we can find transitive models for arbitrarily large finite
  subsets of $T$, and this will allow the rest of the results to carry
  through.  The details are discussed at greater length in Kunen, chs. 4
  and 6.}  And the L\"{o}wenheim-Skolem theorem guarantees that in
addition, we can find such an $M$ that is countable.

Thus, to show that the consistency of ZFC implies that it can't prove
the continuum hypothesis, it will suffice to show that if we have a
countable transitive model $M\models\mathrm{ZFC}$ then we can construct
a countable transitive model
$N\models\mathrm{ZFC}+|\mathbb{R}|=\aleph_\alpha$.  The method used
below will also be able to show that a variety of other hypotheses are
undecidable from ZFC, by showing that models exist of ZFC together with
their negations.

Note that these models make a lot of strange claims.  For instance, if
$M$ is a countable transitive model of ZFC, then we know that
``$\mathbb{R}$ is uncountable''$^M$ is true.  But this seems strange,
because we said the model is countable.  It turns out that what is going
on is that there is a set in $M$ that $M$ ``thinks'' is the set of all
real numbers, but since $M$ only has access to countably many sets, it
only knows about countably many real numbers.  So this set (which I will
call ``$\mathbb{R}^M$'') is in fact countable.  However, since the
sentence ``$x$ is countable'' is written as
$$\exists f((f\colon x\rightarrow \omega)\textrm{ and $f$ is a
  bijection})$$
we see that since $\mathbb{R}^M$ is actually countable, then there is
some $f$ that provides the bijection, but this $f$ is not an element in
$M$.  Thus, $M$ ``gets the cardinality wrong''.  This means that in
general, when talking about cardinalities in a model (or any other
concept that involves implicit quantification, like this), we will always
have to index them with a superscript, as in $\aleph_\alpha^M$.
However, note that if $M\subseteq N$, and $M\models |x|=|y|$ then
$N\models |x|=|y|$ (assuming we are actually talking about the specific
sets $x$ and $y$, rather than using model-relative names like
``$\mathbb{R}$'', or ``$\aleph_\alpha$''), because $N$ has all the
bijections that $M$ does.  So larger models can ``collapse'' certain
cardinalities, but they can't insert new ones between old ones.

It may seem problematic to base our consistency claims about
cardinalities on results about these models that can get cardinalities
so wrong.  But G\"{o}del's results show that if some model satisfies
``$|\mathbb{R}|=\aleph_{17}$'', for instance, then this statement is at
least consistent, even though the relativized statement talks about
different sets $\mathbb{R}^M$ and $\aleph_{17}^M$ rather than the real
ones.  Thus, the strategy described above will be the relevant one - we
just need to be careful when naming infinite cardinalities.

\subsection{$M[G]$}

The more specific strategy will be to start with a countable transitive
model $M\models\mathrm{ZFC}$, find some set $G$ that is not in $M$, and
construct the smallest model $M[G]$ containing every element of $M$ as
well as $G$.  (We will show in section \ref{names} that such a smallest
model exists and is unique.)  To specify $G$ so that we have proper
control over $M[G]$, we will let $P\in M$ be a separative, atomless
partial order (this property is absolute - it can be expressed entirely
in terms of quantifiers ranging over elements of $P$ rather than the
whole universe, so $P$ will actually \emph{be} a separative, atomless
partial order in any model that contains it) and $B$ be the
corresponding complete boolean algebra as described in section
\ref{pos}.  Then $G$ will be a ``generic'' ultrafilter over $B$.  That
is, it will be an ultrafilter such that for every $D\in M$, if $D$ is a
dense subset of $B\setminus\{0\}$, then $G\cap
D\not=\emptyset$.\footnote{This important use of dense subsets is why we
  switch back and forth between partial orders and complete boolean
  algebras.  They are equivalent notions as far as genericity goes, but
  the partial order is easier to describe, and the boolean algebra makes
  some of the technical machinery work more easily.}

Since $M$ is a countable model, we can always guarantee that such a $G$
exists.  This is because $M$ only contains countably many dense subsets
of $B$, and we can number these $D_0,D_1,D_2,\dots$.\footnote{Of course,
  this numbering can only be done externally to the model - inside the
  model, the collection of dense subsets is almost always going to be
  uncountable.}  Since each set is dense, we can find $p_0\in D_0$, and
then $p_1\in D_1$ such that $p_1\leq p_0$, and then $p_2\in D_2$ such
that $p_2\leq p_1$, and so on.  Using the (countable) axiom of choice,
we can fix such a sequence $p_0,p_1,p_2,\dots$, and then let $G=\{x\in
B|\exists i(x\geq p_i)\}$. Since this set contains each $p_i$, we see
that it intersects every $D_i$.  By definition, we can easily see that
$G$ is closed upwards under the ordering.  To see that it is a filter,
we note that if $x,y\in G$, then $x\geq p_i$ and $y\geq p_j$ for some
$i,j$.  Without loss of generality, assume $i\leq j$, so that $x,y\geq
p_j$.  But then $x\wedge y\geq p_j$, so $(x\wedge y)\in G$ and $G$ is a
filter.  We also note that the set $D_p=\{x|x\leq p\text{ or }x\leq\neg
p\}$ is a dense subset of $B$, and $D_p\in M$ since it is definable in
such a simple way.\footnote{Specifically, we just use the axiom schema
  of separation once from $P$, using $\phi(x)$ as``$x\leq p$ or
  $x\leq\neg p$''.  Future such constructions may require more of the
  axioms of ZFC.}  Therefore, $G\cap D_p$ is non-empty, so either $p$ or
$\neg p$ is in $G$, so it is an ultrafilter, as desired.

By choosing $P$ properly, we can make $M[G]$ have certain properties.
For making CH false, we will let $P$ be as follows.  Let $\kappa$ be
some set such that $M\models(|\kappa|=\aleph_\alpha)$,\footnote{If we
  identify cardinals with appropriate ordinals, then we can just let
  $\kappa$ actually \emph{be} $\aleph_\alpha^M$.} for some $\alpha>1$.
Let $P$ be the set of finite partial functions from $\kappa\times\omega$
to $\{0,1\}$.  That is, each element of $P$ specifies finitely many
values for a $\kappa$ by $\omega$ array of zeros and ones.  If $p$ and
$q$ are two elements of $P$, then we will say that $p\leq q$ iff
$q\subseteq p$, so that $p$ is a function extending the function that is
$q$.  Note that this ordering goes in the \emph{opposite} direction from
what one might expect.  It is easy to see that this ordering is
atomless, because any finite partial function can be extended to another
one just by adding one more zero or one to the array.  And to see that
it's separative, note that $p\perp q$ iff there is some location in the
array that one assigns to $0$ and the other assigns to $1$. So if
$x\not\geq y$, so $x$ is not an extension of $y$, then we can find some
value that $x$ assigns and $y$ doesn't, and produce $z$ by switching
that value from $0$ to $1$ or vice versa, and letting $z$ agree with $y$
everywhere else.  Then $x\perp z$ and $z<y$, as required for separativity.

Now, notice the following properties of $G$.  If $p,q\in P\cap G$, then
since $p\wedge q\in G$ and $P$ is dense in $B$, we see that $p\not\perp
q$.  So we see that any two partial functions in $P\cap G$ are
compatible, so their union specifies some (possibly partial) function
from $\kappa\times\omega$ to $\{0,1\}$.  Since $P$ is a dense subset of
$B$, we see that every dense subset of $P$ is a dense subset of $B$ as
well.  Thus, $P\cap G$ intersects every dense subset of $P$.  For any
$x\in\kappa$ and $n\in\omega$, we see that $D_{x,n}=\{p|p(x,n)\text{ is
  defined}\}$ is a dense subset of $P$, because we can always extend any
finite partial function by a single value to get another finite partial
function.  Since $P\cap G$ intersects each of these sets, we see that
the function specified by $P\cap G$ is in fact a total function
$f_G\colon\kappa\times\omega\rightarrow\{0,1\}$.  We can thus define the
functions $G_x\colon\omega\rightarrow\{0,1\}$ by letting
$G_x(n)=f_G(x,n)$.  If $x,y\in\kappa$, then $D_{x,y}=\{p|\exists
n(p(x,n)\not=p(y,n))\}$ is dense as well, because any finite partial
function can be extended by finding some $n$ for which neither $p(x,n)$
nor $p(y,n)$ is defined, and letting one value be $1$ and the other
value be $0$. But then, since $P\cap G$ intersects each of these sets,
we see that $G_x$ and $G_y$ must in fact be different functions from
$\omega$ to $\{0,1\}$.\footnote{We can also see that each of these
  functions must be distinct from any function $M$ knows about.  Let
  $F\in M$ be such that $F\colon\omega\rightarrow\{0,1\}$.  Then let
  $D_{F,x}=\{p|\exists n(p(x,n)\not= F(x)\}$.  Each such set is dense,
  and since $P\cap G$ intersects this set, that means that $G_x\not=F$,
  so each of these functions is distinct from each function that $M$
  had.}  Therefore, since $M[G]$ is a model containing $G$, it can
define each of these functions.  Each such function corresponds to a
subset of $\omega$, so the powerset of $\omega$ must be at least as
large as $\kappa$.  So $M[G]\models(|\mathbb{R}|\geq|\kappa|)$.  Since
$M$ had a bijection between $\kappa$ and $\aleph_\alpha^M$, we see that
$M[G]$ does as well, so $M[G]\models(|\kappa|=|\aleph_\alpha^M|)$.  So
if we can just show that $\aleph_\alpha^M=\aleph_\alpha^{M[G]}$, then we
will have achieved our goal.  I will do this in section \ref{forcing},
but first I will have to go into more detail about the construction of
$M[G]$.

\section{Names and the construction of $M[G]$}\label{names}

Because $G$ is a generic ultrafilter over some boolean algebra $B$ that
$M$ knows about, it turns out that quite a lot about $M[G]$ will be
specifiable just from information about $M$.  In particular, $M$ will
have a ``name'' for every element of $M[G]$, and thus will be able to say
everything that $M[G]$ will be able to.  However, though each such
sentence will have a truth value of either $1$ or $0$ in $M[G]$, $M$
will only be able to specify a ``truth-value'' from $B$.  Sentences
receiving value $0$ or $1$ are known by $M$ to be determinately false or
true, respectively, in $M[G]$.  For all others, only the logical
relations between them are specified, and not the actual truth-values.
But recall from section \ref{BAs} that an ultrafilter in $B$ can be seen
as a way of specifying which elements are to count as ``true''.  It
turns out that if we let $\tv{\phi}=p$ be the ``truth-value'' that $M$
assigns to $p$, then $M[G]\models\phi$ iff $p\in G$.  (This will be
discussed somewhat more in section \ref{forcing}.)  Since there are
many different generic ultrafilters $G$ over $P$, each will give rise to
a slightly different model $M[G]$, and it makes sense that $M$ will have
no way of telling these models apart, since all are equally generic.
But the fact that $M$ can at least say something about the truth-values
of sentences in these models will give enough control to let us show
that $\aleph_\alpha^M=\aleph_\alpha^{M[G]}$ in the case at hand, and to
show other relevant facts in other applications of forcing.

So a $P$-name $n\in M$ will be a set of ordered pairs $(m,p)$, such that
$m$ is a $P$-name, and $p\in P$.  This definition looks circular, since
it defines a $P$-name in terms of $P$-names, but in fact it is
inductive.\footnote{In particular, the circularity is not vicious
  because of the Axiom of Foundation, which ensures that any chain
  $x_0\ni x_1\ni x_2\ni\dots$ eventually terminates after finitely many
  steps.  To know if any set is a $P$-name, it suffices to know which of
  its elements, elements of its elements, and so on, are themselves
  $P$-names.

The universe of sets is constructed by letting $V_0=\emptyset$,
$V_{\alpha+1}=\mathcal{P}(V_\alpha)$, and
$V_\lambda=\bigcup\{V_\beta\colon \beta<\lambda\}$ - every set appears
in some $V_\alpha$ (and all successive ones).  Similarly, the class of
names can be constructed by letting $N_0=\emptyset$, letting
$N_{\alpha+1}$ be the set of all sets of pairs of elements, one from
$N_\alpha$ and one from $B$, and $N_\lambda=\bigcup\{N_\beta\colon
\beta<\lambda\}$.}  I like to think of such a name as specifying a set
of names of potential elements, together with a ``probability'' that each
potential element is actually in the set named.  Then, once the class of
names has been specified, we can specify how to interpret the names as
particular sets.  Given an ultrafilter $G$ over $B$ (and thus a filter
$P\cap G$ over $P$), we will interpret name $n$ as the set
$n^G=\{m^G|(m,p)\in n\text{ and }p\in G\}$.\footnote{Again, because of
  the Axiom of Foundation, we see that this circularity is non-vicious -
  to interpret a name, we just need to be able to interpret all the
  names inside of it, and this process eventually terminates.  Once we
  have interpreted all the names in $N_\alpha$, this rule tells us how
  to interpret the names in $N_{\alpha+1}$, and similarly at limit
  stages.}  Then, once we have fixed the generic ultrafilter $G$, $M[G]$
will be just the set $\{n^G|n\in M\text{ is a $P$-name}\}$.

As an illustration of what the names are like, I will show that
$M\subseteq M[G]$ and that $G\in M[G]$.  For the former, I will
associate with each element $x\in M$ a name $\check x$ such that ${\check
x}^G=x$ for any $G$ whatsoever.  Let $\top$ be the greatest element of
$B$\footnote{In characterizing $B$ as a ring, $\top$ was $1$, but I have
  used a different symbol here to distinguish $\top$ from the natural
  number $1$.  When considering $P$ given above as the set of finite
  partial functions ordered under reverse inclusion, $\top$ is the empty
  function.} and then let $\check x$ be the set $\{(\check y,\top)|y\in
x\}$. Thus, $\check\emptyset=\emptyset$, and if $1=\{\emptyset\}$ then
$\check 1=\{(\check\emptyset,\top)\}=\{(\emptyset,\top)\}$, and if
$2=\{0,1\}$ then $\check 2=\{(\check\emptyset,\top),(\check
1,\top)\}=\{(\emptyset,\top),(\{(\emptyset,\top)\},\top)\}$.  By a
simple induction, one can show that ${\check x}^G=x$.  Because $\top\in G$
for any ultrafilter $G$, we see that ${\check x}^G=\{{\check y}^G|y\in x\}$,
but by our induction assumption, we see that ${\check y}^G=y$, so
${\check x}^G=\{y|y\in x\}=x$.  Thus, every element of $M$ has a name,
so each is in $M[G]$.

To show that $G\in M[G]$, consider the name $\dot X=\{(\check p,p)|p\in
P\}$.  We see that ${\dot X}^G=\{{\check p}^G|(\check p,p)\in \dot
X\text{ and }p\in G\}$.  But ${\check p}^G=p$, and $(\check p,p)\in \dot
X$ iff $p\in P$, so we see that ${\dot X}^G=\{p|p\in P\text{ and }p\in
G\}=P\cap G$.  But it is easy to reconstruct $G$ from $P\cap G$, so
$G\in M[G]$ as well, as desired.\footnote{The construction works as
  follows: $G=\{x\in B|\exists p(p\in P\cap G\text{ and }x\geq p)\}$.
  This can also be done easily within $M[G]$, once we show that it
  satisfies all the axioms of ZFC.}

Now that I have shown that $M\subseteq M[G]$ and $G\in M[G]$, it just
remains to be shown that $M[G]$ is a countable transitive model of ZFC,
as claimed.\footnote{At several points I have also claimed that $M[G]$
  is ``the smallest'' model of ZFC extending $M$ that contains $G$.  But
  any model extending $M$ contains all the $P$-names, and a model
  containing $G$ ``knows how'' to interpret $P$-names (because the
  details of the interpretation can be carried out in any model
  satisfying replacement, separation, and foundation), so it must
  contain all their interpretations, and thus it must contain all of
  $M[G]$.  So I really only need to show that $M[G]$ is in fact a model
  of ZFC, and then it will be clear that it is the smallest such model.}
This proof is not difficult, but it is somewhat tedious.  To get the
idea, note that since the interpretation of any name is a set of
interpretations of names, we see that $M[G]$ is a transitive set, as
described above.  For any transitive set, it is clear that
Extensionality holds in $M[G]$. Because $\check\omega$ is a name, it is
clear that Infinity holds in $M[G]$.  To see that Pairing holds, note
that if $x$ and $y$ are elements of $M[G]$, then there must be names
$\dot x$ and $\dot y$ in $M$ for them.  But then $\dot z=\{(\dot
x,\top),(\dot y,\top)\}$ is a name that denotes their pair.  Separation,
Union, and Powerset work approximately the same way, and Replacement
needs only a slight modification of this technique. The only real
difficulty is in showing that both Foundation and Choice hold in $M[G]$,
but the proof is uninteresting, so I omit it here.\footnote{All these
  proofs can be found in much more detail in any set theory text that
  discusses forcing.}

\section{Forcing}\label{forcing}

Now that we have names for all the elements of $M[G]$, we can use the
relations $\in$ and $=$, together with logical connectives and
quantifiers, to make statements about $M[G]$ in a language expressible
entirely inside $M$.  We will call this language the ``forcing
language'' over $B$.  For any such sentence $\phi$ that contains no
names (like ``$\forall x\forall y(\exists z(\forall w(w\in
z\leftrightarrow(w=x\text{ or }w=y))))$''), we define $M[G]\models\phi$
as before.  But for sentences that contain names, we will need to first
interpret these names relative to $G$.  For instance, if $\dot a,\dot
b,\dot c$ are names in $M$, then
$$\phi=\text{``}\exists x(\dot a\in x\text{ and }\forall y(\dot b\in
y\rightarrow(x\in y\text{ or }y\in \dot c)))\text{''}$$
is a sentence in the forcing language.  In this case, we will say that
$M[G]\models\phi$ just in case
$$\exists (x\in M[G])({\dot a}^G\in x\text{ and }\forall (y\in
M[G])({\dot b}^G\in y\rightarrow(x\in y\text{ or }y\in {\dot c}^G)))$$
With this extension of the satisfaction relation, we can then say that
$p\forces\phi$ (in words, ``$p$ forces $\phi$'') if $M[G]\models\phi$
for every generic ultrafilter $G$ such that $p\in G$.  That is, we will
say that $p\forces\phi$ just in case knowing that $p\in G$ is sufficient
to guarantee that $M[G]\models\phi$.  So sticking $p$ into $G$ is enough
to force $\phi$ to be true in whatever resulting extension $M[G]$ we end
up with.

The amazing and difficult thing to prove is that for every sentence
$\phi$, we can find $\tv{\phi}\in B$ such that $p\forces\phi$ iff
$p\leq\tv{\phi}$.  In addition, it will be the case that $\tv{\text{``It
    is not the case that $\phi$''}}=\neg\tv{\phi}$, and
$\tv{\text{``$\phi$ and $\psi$''}}=\tv{\phi}\wedge\tv{\psi}$, and in
general, logical operations on sentences of the language will correspond
to the operations in the boolean algebra.  This association justifies
calling $\tv{\phi}$ the ``truth-value'' of $\phi$.  The ultrafilter $G$
is then seen as specifying which of these truth-values should be
interpreted as actually being true, for the model $M[G]$ in question.

For every sentence $\phi$ in the forcing language, it will be the case
that $M[G]\models\phi$ iff there is some $p\in G$ such that
$p\forces\phi$.  This fact is what gives $M$ such control over $M[G]$
and allows us to show that $M[G]$ doesn't change too many cardinalities
relative to $M$.  The complete proof of this result is uninteresting and
very long, so I will just refer the reader to Kunen, ch. 6; Jech,
ch. 12; or Bell, ch. 1. \footnote{To start out, note that if $(m,p)\in
  n$, where $m$ and $n$ are names, then $p\forces m\in n$.  However, if
  $(m,p)\in n$ and $(m',p')\in n$, then $m\in n$ might be true even
  though $p\not\in G$ - the reason is because we might have some $q\in
  G$ such that $q\leq p'$ and $q\forces m'=m$.

  To deal with atomic sentences of the form $n=m$, note that if
  $p\forces(\forall x(x\in n\rightarrow x\in m))$ and $q\forces(\forall
  x(x\in m\rightarrow x\in n))$, then $p\wedge q\forces n=m$. With these
  facts (and some more trickery), we can do a simultaneous induction to
  characterize what it takes for $p\forces\phi$ for any atomic sentence
  $\phi$.  (This definition can be found in any of the references.)

  Once we have defined the forcing relation for atomic sentences, a
  further triple induction allows us to characterize the relation for
  all sentences in the forcing language.  The important thing to note is
  that for dealing with quantifiers, we need to use the fact that $B$ is
  a \emph{complete} boolean algebra (in $M$) and not just any boolean
  algebra.  To prove the further facts about truth-values $\tv{\phi}$ we
  need to set up a similar triple induction and show that things work
  out exactly the same way for this induction as for the forcing
  relation.}

Now that we have defined the forcing relation, and know that for every
sentence in the forcing language, if $M[G]\models\phi$ then there is
some $p\in G$ such that $p\forces\phi$, then we can establish further
facts about $M[G]$.  In the example described above, we have shown that
$M[G]\models(|\mathbb{R}|\geq|\aleph_\alpha^M|)$, and we want to show
that $M[G]\not\models\mathrm{CH}$.  As I said before, it will suffice to
show that $\aleph_\alpha^M=\aleph_\alpha^{M[G]}$.  We already know that
$\aleph_\alpha^M\leq\aleph_\alpha^{M[G]}$, because $M[G]$ can only have
more bijections between sets than $M$, not less.  But if
$\aleph_\alpha^M<\aleph_\alpha^{M[G]}$, then there must be some $\beta$
such that $M[G]$ has a bijection between $\aleph_\beta^M$ and
$\aleph_{\beta+1}^M$.  Using the forcing relation, we will be able to
show that this is impossible.

So let $\kappa$ and $\kappa^+$ be some sets in $M$ that $M$ assigns
cardinalities $\aleph_\beta^M$ and $\aleph_{\beta+1}^M$.  Because $M[G]$
extends $M$, these sets will be in $M[G]$ as well.  By the
countable Axiom of Choice, we know that a union of $\aleph_\beta$ many
countable sets is itself of size $\aleph_\beta$, and not
$\aleph_{\beta+1}$, so if we can find some function $F\in M$ that
assigns each element of $\kappa$ to a countable subset of $\kappa^+$,
then the union of these subsets must leave out some element of
$\kappa^+$.  Now, let us assume for the sake of contradiction that there
is some function $f\in M[G]$ (with name $\dot f$) such that
$$M[G]\models\text{``$\dot
  f\colon\check\kappa\rightarrow\check{\kappa^+}$ is a bijection.''}$$
By the facts about forcing, we see that there is $p\in G$ such that
$$p\forces\text{``$\dot
  f\colon\check\kappa\rightarrow\check{\kappa^+}$ is a bijection.''}$$
Now, I will define $F_p\in M$ as a function from $\kappa$ to subsets of
$\kappa^+$ as follows.  For any $y\in\kappa^+$, we will have $y\in
F_p(x)$ iff there is $p_{x,y}\leq p$ such that
$$p_{x,y}\forces\text{``}\dot f(\check x)=\check y\text{''}.$$
Because $P$ is dense in $B$, we can assume that each of these $p_{x,y}$
is an element of $P\subseteq B$.  But we can see that $p_{x,y}$ must be
incompatible with $p_{x,y'}$ if $y'\not=y$, since these conditions both
force $\dot f$ to be a bijection, but force it to take different values
at the same point.  But any collection of incompatible elements of $P$
must be countable in $M$.\footnote{This is because elements of $P$ are
  finite partial functions, and they are incompatible iff they assign
  different values to the same input.  There is a combinatorial result
  called the ``$\Delta$-system lemma'' (see Kunen or Jech) that states
  that for any uncountable collection of finite sets, there is an
  uncountable subcollection such that any pair of them has the same
  intersection, and we call this intersection the ``root'' of the
  $\Delta$-system.  But if any two of these functions are incompatible,
  they must disagree somewhere in the root (because if they disagree
  only elsewhere, then they are in fact compatible).  But there are only
  finitely many sets of values that can be assigned on the root, so the
  original collection of pairwise incompatible elements must have been
  countable.}  Thus, each of these sets $F_p(x)$ must be countable in
$M$.  Therefore, their union cannot be all of $\kappa^+$.  But because
$f\in M[G]$ is in fact a bijection between $\kappa$ and $\kappa^+$, we
see that for every $y\in\kappa^+$, there is in fact some $x\in\kappa$
and $p_{x,y}\in G$ such that $p_{x,y}\leq p$ and
$p_{x,y}\forces\text{``}\dot f(\check x)=\check y\text{''}$.  Thus,
every element of $\kappa^+$ is in some $F_p(x)$.  This is a
contradiction, so our original assumption (that there was $f\in M[G]$
that was a bijection between $\kappa$ and $\kappa^+$) must have been
false.  Therefore, $M[G]$ has all the same cardinalities as $M$, as
required, so our earlier construction in fact shows that
$$M[G]\models|\mathbb{R}|\geq\aleph_\alpha^{M[G]}$$
which means that $M[G]$ falsifies the Continuum Hypothesis as badly as
we want.  Therefore, the negation of CH is consistent with ZFC, QED.

\section{Other applications}

\subsection{Making CH true}

We can use a different partial order $P\in M$ to show that we can make
models where CH is true, in addition to models where CH is false.  In
this case, we let $P$ be the set of \emph{countable} (in $M$) partial functions
from $\aleph_1^M$ to $\mathbb{R}^M$, and again we say that $p\leq q$ iff
$p$ is a function extending $q$.  Again, we see that this order is
atomless, because any countable partial function can be extended, and it
is separative, because if $p<q$, then we can find some value that $p$
assigns and $q$ doesn't, and change it to get an incompatible countable
partial function extending $q$.

Now, just as before, since $G$ is a filter, we see that it defines a
partial function from $\aleph_1^M$ to $\mathbb{R}^M$.  For any
$x\in\aleph_1^M$, we can define $D_x=\{p|p(x)\text{ is defined}\}$ to
see that this function is in fact total.  Defining $D_r=\{p|\exists
x(p(x)=r)\}$ we can see that this function is a surjection.  Thus,
$M[G]\models|\mathbb{R}^M|\leq|\aleph_1^M|$.  We already know that
$\aleph_1^M\leq\aleph_1^{M[G]}$.  So if we can show that
$\mathbb{R}^M=\mathbb{R}^{M[G]}$, then we will see that
$M[G]\models\text{CH}$.  To do this, we must show that if $f\in M[G]$ is
a function $f\colon\omega\rightarrow\{0,1\}$, then $f\in M$.

So let $\dot f$ be a name in $M$ for $f$.  Choose $p\in G$ such that
$p\forces\dot f\colon\check\omega\rightarrow\{\check 0,\check 1\}$.  Now
find a $p_0\in G$ forcing the value of $\dot f$ at $0$ such that
$p_0\leq p$. That is, either $p_0\forces\dot f(\check 0)=\check 0$ or
$p_0\forces\dot f(\check 0)=\check 1$.  Similarly, find $p_1\leq p_0$ in
$G$ forcing the value of $\dot f$ at $1$, and $p_2\leq p_1$ forcing the
value of $\dot f$ at $2$, and so on.  Then we get a descending sequence
of elements of $P$ that force $\dot f$ to be a function from $\omega$ to
$\{0,1\}$ and together force all of its values.  Since these are a
descending sequence of countable functions, their union is itself a
countable function, so call this element of the order $q$.  (This is why
we had to use countable partial functions, rather than finite partial
functions as before.)  Now, since $q\in M$ and $\dot f\in M$, and $q$
forces every value of $\dot f$, we see that $M$ ``knows how'' to read
off all the values of $\dot f$ from $q$.  That is, a single element of
the partial order was sufficient to specify it, rather than needing the
entire ultrafilter.  So the function $f$ must be in $M$ as well as in
$M[G]$.  Therefore, $\mathbb{R}^{M[G]}=\mathbb{R}^M$, as required, so
$M[G]\models\text{CH}$.

\end{document}